\documentclass[12pt]{article} 
\begin{document}

\def\C{{\bf C}}
\def\Q{\hbox{\rm\ l\hskip -1.7truemm Q}}
\def\R{\hbox{\rm\ l\hskip -0.7truemm R}}
\def\P{\hbox{\rm\ l\hskip -0.7truemm P}}
\def\N{\hbox{\rm\bf N}}
\def\Z{\hbox{\rm\bf Z}}
\def\ord{{\rm ord} }
\def\mult{{\rm mult}}
\def\endProof {\rule[-0.5mm]{1.2ex}{1.2ex}}

\author{ Nguyen Van Chau
\thanks {  Supported in part by the National Basic Program on Natural Science, 
Vietnam.}  }

\title{Non-proper value set and the Jacobian condition} 
\date {}
\maketitle

\begin{abstract}
The non-proper value set of a nonsingular polynomial map from $\C^2$ into 
itself, if non-empty, must be a curve with one point at infinity.

2000{ \it Mathematical Subject Classification:} 14 H07,  14R15.
  
{\it Key words and phrases:} non-proper value set, Jacobian condition.
 
\end{abstract}

\medskip

\noindent{\bf 1.} Let $f=(P,Q):\C^2_{(x,y)}\longrightarrow \C^2_{(u,v)}$ be a 
dominant polynomial map, $P,Q\in \C [x,y]$ and denote $J(P,Q):=P_xQ_y-
P_yQ_x$. Recall that the so-called {\it non-proper value set } $A_f$ of $f$ is  
the set  consists of all point $a\in \C^2$ such that the inverse $f^{-1}(K)$ is not 
compact for compact neighbourhoods $K\subset \C^2$ of $a$. This set $A_f$, if 
non-empty, must be a plane curve that each of it's irreducible components can be 
parameterized by a non-constant polynomial map from $\C$ into $\C^2$ ( See 
[J]). The mysterious Jacobian conjecture (See [BMW] and [E]), posed first by 
Keller in 1939 and still open,  asserts that a polynomial map $f=(P,Q)$ of $\C^2$ 
with $J(P,Q)\equiv const. \neq 0$ must have a polynomial inverse. This 
conjecture can be reduced to prove that the non-proper value set $A_f$ is empty. 
In any way one may think that in a counterexample to the Jacobian conjecture, 
if exists, the non-proper value set must has a very special configure. Such 
knowledges  may be useful in pursuit of  this conjecture. In [C]  it was observed 
that   the irreducible components of $A_f$ in such a counterexample  can be 
parameterized by polynomial maps $\xi\mapsto (p(\xi),q(\xi))$ with $\deg p/\deg 
q=\deg P/ \deg Q$. In this  paper we would like to notice that the non-proper 
value set of a nonsingular polynomial map from $\C^2$ into itself, if non-empty, 
must be a curve with one point at infinity.

\medskip

\noindent{\bf Theorem 1. }{\it Suppose $f=(P,Q)$ is a polynomial map of $\C^2$ 
with $J(P,Q)\equiv const. \neq 0$, $\deg P=Kd$ and $\deg Q=Ke$, $\gcd (d,e)=1
$, $P$ and $Q$ are monic in $y$,
$$
\begin{array}{c}
P(x, y)=Ay^{Kd}+\dots +a_1(x)y+a_0(x), \ A\neq 0\cr
Q(x, y)=By^{Ke}+\dots +b_1(x)y+b_0(x), \ B\neq 0.
\end{array}
\eqno (1)
$$
If  the non-proper value set  $A_f$ is not empty, then every irreducible 
component of $A_f$ can be parameterized by polynomial maps of the form 
$$\xi\mapsto (A\xi^{md}+\mbox { \rm lower terms in } \xi, B\xi^{me}+\mbox { 
\rm lower terms in } \xi), \ m \in {\bf N}.\eqno(2)$$
}

\medskip
By definition $A_f$ is the set of all value $a\in \C^2$ such that the number of 
solutions counted with multiplicities of the equation $f(x,y)=a$ is different from 
those for generic values in $\C^2$. Then, considering the components $P(x,y)$ 
and $Q(x,y)$ as elements of $\C[x][y]$ we can define the resultant 
$$Res_y(P-u,Q-v)=R_0(u,v)x^N+\dots+R_N(u,v),\eqno (3)$$
where $ R_i\in \C[u,v],\ R_0\neq 0.$ From the basic properties of the resultant 
function we know that $N$ is the geometric degree of $f$ and $A_f=\{(u,v)\in 
\C^2: R_0(u,v)=0\}$. Note that a curve given by a polynomial parameter of the 
form (1) can be defined by a polynomial of the form  $(A^eu^e-B^d v^d)^m+ 
\sum_{0\leq id+je <mde} c_{ij}u^iv^j$ and its branch at infinity has a Newton-
Puiseux series of the form $u=cv^{d\over e}+\mbox { \rm lower terms in } v$, 
where $c$ is a   $d-$ radicals of $B^d/A^e$ . Thus, Theorem 1 leads  to  

\medskip
\noindent{\bf Corollary 1. }{\it Let $f$ be as in Theorem 1. Then
$$
R_0(u,v)=C(A^eu^e-B^dv^d)^M+ \sum_{0\leq id+je <Mde} c_{ij}u^iv^j\eqno 
(4)
$$
with $\ 0\neq C\in \C$ and $\ M \geq 0.$ 
}

\medskip 

\medskip
\noindent{\bf Corollary 2. }{\it Let $f$ be as in Theorem 1. If $A_f \neq 
\emptyset$, then $A_f$ is a curve with one point at infinity and the  irreducible 
branches at infinity of $A_f$  have  Newton-Puiseux series of the form 
$$u=cv^{d\over e}+\mbox { \rm lower terms in } v$$ 
with  coefficients  $c$ to be  $d-$ radicals of $B^d/A^e$. }

\medskip 
As seen later, the monic representation in (1) of $P$ and $Q$ is only used to 
visualize the coefficient $B^d/A^e$. In fact, when $A_f\neq \emptyset$ the 
numbers $d$, $e$, $B^d/A^e$ and the polynomial $R_0(u,v)$ are invariant  of 
$f$ under right actions of automorphisms of $\C^2$, since  the set $A_f$ does 
not depend on the coordinate $(x,y)$. Furthermore, the coefficient $B^d/A^e$ is 
uniquely determined in the relation
$$P_+^e(x,y)=(B^d/A^e)Q_+^d(x,y),$$
which is dominated by the Jacobian condition when $\deg P >1$ and $\deg Q >1
$. Here,  $P_+$ and $Q_+$ are leading homogenous components of $P$ and 
$Q$, respectively.

\medskip
Theorem 1 will be proved  in sections 2 - 5  by an elementary way using Newton-
Puiseux expansions and the Newton Theorem. It is worth to determine the form 
of  $R_0(u,v)$ by examining directly the resultant function $Res_y(P-u,Q-v)$.

\medskip

\noindent{\bf 2. Dicritical series of $f$.} In order to prove Theorem 1, we need to 
setup our hypothesis. From now on, $f=(P,Q):\C^2_{(x,y)} \longrightarrow \C^2_
{(u,v)}$ is a given  polynomial map with $J(P,Q)\equiv const. \neq 0$, $\deg P 
=Kd >0 $ and $\deg Q=Ke >0 $, $\gcd (d,e)=1$. The Jacobian condition will be 
used really in Lemma 3 and the proof of Theorem 1.  Since $A_f$ does not 
depend on the coordinate $(x,y)$, to examine it we can assume that  $P$ and 
$Q$  are monic in $y$, 
$$
\begin{array}{c}
P(x, y)=Ay^{Kd}+\dots +a_1(x)y+a_0(x), \ A\neq 0\cr
Q(x, y)=By^{Ke}+\dots +b_1(x)y+b_0(x), \ B\neq 0.
\end{array}
\eqno (5)
$$
With this representation the Newton-Puiseux roots at infinity $y(x)$ of each 
equations $P(x,y)=0$ and $Q(x,y)=0$ are fractional power series of the form
$$y(x)=\sum_{k=0}^\infty c_kx^{1-{k\over m}}, \ m\in \N, \ \gcd\{ k:c_k\neq 0
\}=1,$$
for which the map $\tau \mapsto (\tau^m,y(\tau^m))$ is meromorphic and 
injective for $\tau$ larger enough .  In view of the Newton theorem we can 
represent
$$P(x,y)=A\prod_{i=1}^{\deg P}(y-u_i(x)), \
Q(x,y)=B\prod_{j=1}^{\deg Q}(y-v_i(x)), \eqno (6)$$
where $u_i(x)$ and $v_j(x)$ are the Newton-Puiseux roots at infinity of the 
equations $P=0$ and $Q=0$, respectively. We refer the readers to [A] and [BK] 
for the Newton theorem and the Newton-Puiseux roots.

We begin with the description of the non-proper value set  $A_f$ of $f$ via 
Newton-Puiseux expansion. We will work with finite fractional power series 
$\varphi(x,\xi)$ of the form 
$$\varphi (x,\xi)=\sum_{k=1}^{K-1} a_kx^{1-{k\over m}}+\xi x^{1-{K\over 
m}}, m\in \N 
,\ \gcd\{k: a_k \neq 0 \}=1, \eqno (7)$$
where $\xi$ is a  parameter. For convenience, we denote $\mult (\varphi):=m$. A  
such  series $\varphi$ is called  {\it dicritical series } of $f$ if  
$$f(x, \varphi (x, \xi))=f_\varphi (\xi) +{\mbox {\rm lower terms in }} x, \  \deg 
f_\varphi >0.$$ 
The following description of $A_f$ was presented in [C].

\medskip
\noindent{\bf Lemma 1.} ( Lemma 4, [C])
$$ A_f = \bigcup_{\varphi  \ is\ a\ dicritical\ series\ of\  f} f_{\varphi} (\C). $$

\medskip
To see it, note that by definitions  the non-proper value set  $A_f$ consists of all 
values $a\in \C^2$ such that there exists a sequence $\C^2\ni p_i\rightarrow 
\infty$ with $f(p_i)\rightarrow a$. If $\varphi$ is a dicritical series of $f$ of the 
form (7), we  can define the map $\Phi(t,\xi):=(t^{-m},\varphi(t^{-m},\xi))$. 
Then, $\Phi$ sends $\C^*\times\C$ to $\C^2$  and the line $\{ 0\}\times \C $ to 
the line at infinity of $\C\P^2$. The polynomial map $F_\varphi (t, \xi):=f\circ 
\Phi(t, \xi)$  sends the line $\{ 0\}\times \C$ into $A_f \subset\C^2. $ Therefore, 
$f_\varphi(\C)$ is an irreducible component of $A_f $, since  $\deg f_\varphi >0
$. Conversely, if $\ell$ is an irreducible component of $A_f$, one can choose a 
smooth point $(u_0,v_0)$ of $A_f$, $(u_0,v_0)\in\ell$, and  an irreducible 
branch at infinity $\gamma$ of the curve $P=u_0$ (or the curve $Q=v_0$) such 
that the image $f(\gamma)$ is a branch curve  intersecting transversally  $\ell$ 
at $(u_0, v_0).$ Let $u(x)$ be a Newton-Puiseux expansion at infinity of 
$\gamma$. Then, we can construct an unique dicritical  series $\varphi (x,\xi)$ 
such that $u(x)=\varphi (x,\xi_0+\mbox{ lower term in } x) $. For this dicritical 
series $\varphi$ we have  $f_\varphi (\C)=\ell$.

\medskip
\noindent{\bf 4.} {\it Associated sequence of dicritical series.} 
Let  $\varphi$ be a given dicritical series of $f$.  Let us to represent 
$$ \varphi (x, \xi)= \sum_{ k=0}^{ K-1} c_kx^{1-\frac{n_k}{m_k} }+\xi x^{ 1-
\frac{ n_K}{ m_K}},  \eqno (8)$$
where $0\leq \frac{n_0}{m_0} <\frac{n_1}{m_1}< \dots <\frac{n_{K-1}}{m_
{K-1}}< \frac{n_K}{m_K}=\frac{n_\varphi}{ m_\varphi}  $ 
and  $c_i\in \C$ may be the zero,  so that  the sequence of series $ \{\varphi_i \}_
{ i=0,1\dots,K}$ defined by  
$$\varphi_i(x, \xi ):=\sum_{ k=0}^{ i-1} c_kx^{1-\frac{n_k}{m_k}}+\xi x^{1-
\frac{n_i}{m_i}},  i=0, 1,\dots ,  K-
1, \eqno (9)$$ 
and $\varphi_K:=\varphi$ satisfies the following properties: 

S1) $\mult (\varphi_i)=m_i$.

S2) For every $i < K$ at least one of polynomials $p_{\varphi_i}$ and $q_
{\varphi_i}$ 
has a zero point different from the zero.  
 
S3) For every $ \psi (x, \xi)=\varphi_i(x, c_i)+\xi x^{1-\alpha}$,  $ \frac{n_i}
{m_i} <\alpha < \frac{n_{i+1}}{m_{i+1}}, $ each of the polynomials $p_\psi$ 
and $q_\psi$ is either constant or a monomial of $\xi. $ 
 
The representation (8) of $\varphi$ is thus the longest representation such that 
for each index $i$ there is a Newton-Puiseux root $y(x)$ of $P=0$ or $Q=0$ 
such that
$y(x)=\varphi_i(x,c+\mbox{ lower terms in } x)$, $c\neq 0$ if $c_i=0$. This 
representation and  the associated sequence $ \{\varphi_i \}_{ i=0,1,\dots ,K}$ 
are  well defined and unique. Further, $\varphi_0(x,\xi)=\xi x$.

We will use the associated sequence $\{\varphi_i\}$ to determine the form of the 
polynomials $f_\varphi(\xi)$. For simplicity in notations,  in below we shall  use 
lower indeces ``$i$" instead of the lower indeces ``$\varphi_i$".

For each  associated series $\varphi_i$, $i=0.\dots, K$, let us represent
$$
\begin{array}{c}
P(x, \varphi_i (x, \xi))= p_i(\xi)x^\frac{ a_i}{ m_i}+\mbox {\rm lower terms in } 
x \\
Q(x, \varphi_i (x, \xi))=q_i (\xi)x^\frac{ b_i}{ m_i }+\mbox {\rm lower terms 
in } x,\\
\end{array}
\eqno (10)$$
where $ p_i, q_i \in \C[\xi]-\{0\}$,   $a_i,b_i\in \Z$ and $m_i:=\mult (\varphi_i)
$. 

\medskip

Let $\{ u_i(x), i=1,\dots \deg P\}$ and $\{v_j(x),j=1,\dots \deg Q\}$ be the 
collections of the Newton-Puiseux roots of $P=0$ and $Q=0$, respectively. As 
shown in Section 2, by  the Newton theorem the polynomials $P(x,y)$ and $Q
(x,y)$ can be factorized in the form
$$P(x,y)=A\prod_{i=1}^{\deg P}(y-u_i(x)),  Q(x,y)=B\prod_{j=1}^{\deg Q}(y-
v_j(x)).\eqno (11)$$
For each  $i=0.\dots, K$, let us define 

- $S_i:=\{k: 1\leq k\leq \deg P: u_k(x)=\varphi_i(x,a_{ik}+\mbox{ lower terms 
in } x) a_{ik}\in\C\}$;

- $T_i:=\{k: 1\leq k\leq \deg Q: v_k(x)=\varphi_i(x,b_{ik}+\mbox{ lower terms 
in } x), b_{ik}\in\C\}$;

- $S_i^0:=\{k \in S_i: a_{ik}=c_i\}$;

 - $T_i^0:=\{k \in T_i: b_{ik}=c_i\}.$\\
Represent

$$p_i(\xi)=A_i\bar p_i(\xi)(\xi-c_i)^{\# S_i^0}, \bar p_i(\xi):=\prod_{k\in 
S_i\setminus S_i^0}(\xi-a_{ik}),$$
and
$$q_i(\xi)=B_i\bar q_i(\xi)(\xi-c_i)^{\# T_i^0}, \bar q_i(\xi):=\prod_{k\in 
T_i\setminus T_i^0}(\xi-b_{ik}).$$
\medskip
\noindent{\bf Lemma 2. }{\it 

i) $n_0=0$, $m_0=1$ and 
$$A_0=A, \deg p_0=a_0=Kd$$ 
$$B_0=B,  \deg q_0=b_0=Ke.$$

ii) For $i=1,\dots , K$ 
$$A_i=A_{i-1}\bar p_{i-1}(c_{i-1}), \deg p_i=\# S_i=\#S_{i-1}^0$$ 
$${a_i\over m_i}={a_{i-1}\over m_{i-1}} + \# S_{i-1}^0(\frac{n_{i-1}}{m_{i-
1}}-\frac{n_i}{m_i})$$

$$B_i=B_{i-1}\bar q_{i-1}(c_{i-1}), \deg q_i=\# T_i=\#T_{i-1}^0,$$
$${b_i\over m_i}={b_{i-1}\over m_{i-1}} + \# T_{i-1}^0(\frac{n_{i-1}}{m_{i-
1}}-\frac{n_i}{m_i}),$$

}
\medskip
{\it Proof.} Note that $\varphi_0(x,\xi)=\xi x$ and $\varphi_i(x,\xi)=\varphi_{i-
1}(x,c_{i-1})+\xi x^{1-\frac{n_i}{m_i}}$ for $i>0$. Then, substituting 
$y=\varphi_i(x,\xi)$, $i=0,1,\dots , K$,  into the Newton factorizations of $P(x,y)
$ and $Q(x,y)$ in (11) one can easy verify the conclusions. \endProof

\medskip
\noindent{\bf 4. The Jacobian condition.} Let $\varphi $ be a dicritical series of 
$f$ and $\{ \varphi_i\}$ be it's associated series. 
Denote $$J_i(\xi):=a_ip_i(\xi)\dot q_i(\xi)-b_i\dot p_i(\xi) q_i(\xi).$$
The Jacobian  condition will be considered in the following meaning. 

\medskip
\noindent{\bf Lemma 3}: {\it Let $0\leq i <K$.  If $a_i >0$ and $b_i>0$, then
$$ 
J_i(\xi) \equiv \cases{-m_i J(P,Q)&if $a_i+b_i=2m_i-n_i,$ 
\cr 
        0&if $a_i+b_i>2m_i-n_i.$\cr} 
$$
Further, $J_i(\xi)\equiv 0$ if and only if $p_i(\xi)$ and $q_i(\xi)$ have a 
common zero point. In this case   
$$p_i(\xi)^{b_i}=Cq_i(\xi)^{a_i}, \ C\in \C^*.$$
}
 
\medskip
\noindent{\it Proof.} Since $a_i >0$ and $ b_i> 0,$  taking  differentiation of 
$Df(t^{-m_i},\varphi_i (t^{-m_i},\xi ))$,  we have that  
$$m_i J(P,Q)t^{n_i-2m_i-1} +\mbox { higher terms in }t    
=-J_i (\xi) t^{-a_i-b_i-1}+ \mbox { higher terms in }t. 
$$ 
Comparing two sides of it we  can get the first conclusion.
The remains are left to the readers as an elementary exercise. \endProof

\medskip
\noindent{\bf 5. Proof of Theorem 1.}

i) Assume that $A_f\neq \emptyset$. Then,  $A_f$ is a plane curve in $\C^2$. 
Let $\ell$ be an irreducible component of $A_f$. By Lemma 1
there is a dicritical series $\varphi$ of $f$ such that $\ell$ can be parameterized 
by the polynomial map $f_\varphi(\xi)=(p_\varphi(\xi),q_\varphi(\xi))$, i.e. 
$\ell=f_\varphi(\C)$. We will show that
$$
f_\varphi(\xi)=(AC^d_\varphi\xi^{ D_\varphi  d}+\dots,
BC^e_\varphi \xi^{ D_\varphi  e}+\dots  ),\ 
C_\varphi \neq 0,\ D_\varphi \in \N.\eqno (12)$$ 
Then, by changing variable $\xi\mapsto C_\varphi^{-1}\xi$ we get the  desired 
parameterization $\xi\mapsto (A\xi^{D_\varphi d}+\dots,B\xi^{D_\varphi e}
+\dots  ) $ of $\ell$.

\medskip

ii) Consider the associated sequence $\{\varphi_i\}_{i=1}^K$ of $\varphi$.  
Since  $A_f\neq \emptyset$ as assumed, $$\deg P>1, \deg Q>1.$$
Otherwise, $f$ is bijective and $A_f=\emptyset$. Since $\varphi$ is a dicritical 
series of $f$,  without loss of generality we can assume that
$$\deg p_K >0, \ a_K=0 \mbox { and } b_K\leq 0.$$
Then, from the constructing of the sequence $\varphi_i$ it follows that  
$$
\cases{
p_i(c_i)=0 \mbox { and } a_i>0 ,& $ \ i=0,1,\dots ,K-1$ \cr
q_i(c_i)=0&  if   $ b_i >0$}
\eqno (13)$$
This allows us to use the Jacobian condition in the meaning of Lemma 3. Then, 
by   induction using Lemma 2, Lemma 3 and (13) we can obtain without 
difficult
the following.

\medskip
\noindent{\it Assertion:}
{\it For $i=0,1,\dots ,K-1$ we have
$$a_i>0,b_i>0,\eqno (a)$$
$${a_i\over b_i}={\#S_i\over \# T_i}={d\over e}\eqno(b)$$
and$${\# S_i^0\over \#T_i^0}={d\over e},
\bar p_i(\xi)^e=\bar q_i(\xi)^d . \eqno(c)$$}

\medskip

iii) Now, we prove (12). By Lemma 2 (iii) and (b-c) we have

$\frac{b_K}{m_K}$

$=\frac{b_{K-1}}{m_{K-1}}+\#T^0_{K-1}(\frac{n_{K-1}}{m_{K-1}}-\frac
{n_K}{m_K})$

$=\frac{e}{d}[\frac{a_{K-1}}{m_{K-1}}+\#S^0_{K-1}(\frac{n_{K-1}}{m_{K-
1}}\frac{n_K}{m_K})]$

$=\frac{e}{d}\frac{a_K}{m_K}$

$=0,$ 

\noindent as $a_K=0$.
Hence,  $f_\varphi(\xi)=(p_K(\xi),q_K(\xi))$ by definition and (a). Using 
Lemma 2 (ii-iii) to compute the coefficient $A_K$ and $B_K$  we  can get
$$A_K=A(\prod_{k\leq K-1}\bar p_k(c_k))
, B_K=B(\prod_{k\leq K-1}\bar q_k(c_k)).$$
Let $C_\varphi$ be a $d-$radical of $(\prod_{k\leq K-1}\bar p_k(c_k))$ and 
$D_\varphi:=\gcd (\# S^0_{K-1},\# T^0_{K-1})$. Then, by Lemma 2 (ii) and (b-
c) we have
that $A_K=AC_\varphi^d$,  $B_K=BC_\varphi^e$, $\deg p_K=\# S^0_{K-1}
=D_\varphi d$ and $ \deg q_K=\# T^0_{K-1}=D_\varphi e$. Thus,
$$
f_\varphi(\xi)=(AC^d_\varphi\xi^{ D_\varphi  d}+\dots,
BC^e_\varphi \xi^{ D_\varphi  e}+\dots  ).$$ 
\endProof

\medskip

\bigskip

{\it Acknowledgment.} The author wishes to thank Prof. H.V. Ha  for many 
valuable suggestions  and useful discussions.

\bigskip
 
\noindent {\bf References} 

 \noindent [AM] S.S.  Abhyankar and T.T.   Moh,     {\it  Embeddings of the line 
in the 
plane},     J.   Reine Angew.   Math.   276 (1975),     148-166.   

\noindent [A]  S. S. Abhyankar,    {\it Expansion Techniques in Algebraic 
Geometry},    
Tata 
institute of fundamental research,    Tata Institute,    1977.  

\noindent [BCW] H.   Bass,     E.   Connell and D.   Wright,    {\it  The Jacobian 
conjecture: 
reduction of degree and formal expansion of the inverse},     Bull.   Amer.   
Math.   
Soc.   (N.S.) 7 (1982),     287-330.   

\noindent [BK] E.  Brieskorn,    H.  Knorrer,    {\it Ebene algebraische Kurven},    
Birkhauser,    Basel-Boston-Stuttgart 1981.  

\noindent [C] Nguyen Van Chau,  {\it 
Non-zero constant Jacobian polynomial maps of $ C^2. $}
J.  Ann.  Pol.  Math.  71,  No. 3,  287-310 (1999).  

\noindent [E] van den Essen,  Arno,  {\it
Polynomial automorphisms and the Jacobian conjecture}.  (English.  English 
summary) Progress in Mathematics,  190.  Birkhäuser Verlag,  Basel,  2000. 

\noindent [J] Z.  Jelonek, {\it The set of points at which a polynomial map is not 
proper.} Ann.  Pol.  Math,  58 (1993),  259-266.

\bigskip
\noindent{\small Hanoi Institute of Mathematics,
P. O.  Box 631, Boho 10000,  Hanoi,  Vietnam. \\E-mail: 
nvchau@thevinh.ncst.ac.vn
\end{document}